\documentclass[conference]{IEEEtran}
\IEEEoverridecommandlockouts
\usepackage{cite}
\usepackage{amsmath,amssymb,amsfonts}
\usepackage{algorithmic}
\usepackage{graphicx}
\usepackage{textcomp}
\usepackage{xcolor}
\usepackage{algorithm}
\usepackage{algorithmic}
\usepackage{hyperref}
\usepackage{listings}
\usepackage[dvipsnames]{xcolor}

\usepackage{listings}

\lstdefinestyle{mystyle}{
    backgroundcolor=\color{white},   
    commentstyle=\color{gray},
    keywordstyle=\color{blue},
    numberstyle=\tiny\color{gray},
    stringstyle=\color{teal},
    basicstyle=\ttfamily\footnotesize,
    breaklines=true,
    captionpos=b,
    numbers=left,
    numbersep=5pt,
    showstringspaces=false,
    tabsize=2,
}

\def\BibTeX{{\rm B\kern-.05em{\sc i\kern-.025em b}\kern-.08em
    T\kern-.1667em\lower.7ex\hbox{E}\kern-.125emX}}
\begin{document}

\title{Reproducing the $k$-copwin Algorithm: Theory vs. Implementation}

\author{\IEEEauthorblockN{Corresponding Author: Meagan Mann}
\IEEEauthorblockA{
\textit{Queen's University}\\
Kingston, Canada \\
20mmm27@queensu.ca}
\and
\IEEEauthorblockN{Christian Muise}
\IEEEauthorblockA{
\textit{Queen's University}\\
Kingston, Canada \\
christian.muise@queensu.ca}
\and
\IEEEauthorblockN{Erin Meger}
\IEEEauthorblockA{
\textit{Queen's University}\\
Kingston, Canada \\
erin.meger@queensu.ca}

}

\maketitle

\begin{abstract}

Cops and Robbers is a well-studied pursuit-evasion game that provides insights into graph theory and theoretical computing. A central question is determining the minimum number of cops required to capture the robber, known as the cop number. We focus on reproducing an algorithm proposed by Petr, Portier, and Versteegen in 2022, which efficiently determines whether a graph is $k$-copwin. This paper presents a Python implementation of the $k$-copwin algorithm. In this work, we present our implementation in detail, clarify key aspects of the algorithm, and discuss its implications for future practical deployments.

\end{abstract}


\section{Introduction}

\subsection{Background}

The game of Cops and Robbers was introduced independently by Nowakowski and Winkler \cite{NOWAKOWSKI1983235} and Quillot \cite{quilliot1978jeux}. Cops and Robbers is a pursuit-evasion game played on a graph $G$, where $k$-cops try to capture the robber, by having at least one cop occupy the same vertex as the robber. When this process is successful, we say the graph $G$ is $k$-copwin. The cops are first placed anywhere on the graph, followed by the robber, then they each take turns making a move by either remaining on the current vertex or moving to a neighbouring vertex. Multiple cops can occupy the same vertex at any point during the game. We say that $c(G)$ is the cop number of a graph $G$, and is the minimum $k$ such that $G$ is $k$-copwin.

Over the past four decades, Cops and Robbers has evolved into an active research area intersecting structural graph theory, combinatorics, and algorithm design. Early structural results from Aigner and Fromme \cite{fromme1984game} include that planar graphs have cop number at most 3, showing that cop number is related to topological structure. From there, research expanded to explore how cop number relates to graph parameters such as genus \cite{BERARDUCCI1993389}, girth \cite{https://doi.org/10.1002/jgt.22855}, treewidth \cite{Kosowski2015KChordal}, diameter \cite{BonatoHahnWang2007CopDensity}, and an influential open problem called Meyniel's conjecture that states for a connected graph of size n $c(g) = O(\sqrt{n})$ \cite{frankl1987cops,baird2013meynielsconjecturecopnumber}. Connections between cop number and graph minors show that cop number is closely tied to structural decomposition of the graph. Recent work on graph minors and cop number \cite{Forbidminors} shows how forbidden minor structure constrains pursuit dynamics. Similarly, work on induced subgraphs \cite{clow2025copsrobbersgraphspath} shows how cop number behaves under induced substructure, reinforcing that the game captures combinatorial properties rather than just local adjacency relations. More recent work has also looked at extremal constructions such as high-girth \cite{highgirth}, high-minimum-degree and cage graphs \cite{TURCOTTE202174}, where the cop number behavior becomes highly dependent on global structure. These examples often require very large graphs before asymptotic structural patterns emerge. Today, Cops and Robbers remains an active area of research in extremal bounds, structural characteristics, and algorithmic developments.

\subsection{Algorithmic Motivation}
The Cops and Robbers game has served as a foundation for exploring both the problems in pursuit-evasion theory and the algorithmic techniques required to solve them. Clarke and MacGillivray \cite{clarke2012characterizations} provided an algorithm  with complexity $O(n^{2k+2})$, where $k$ denotes the number of cops, and $n$ denotes the number of vertices in the graph. While polynomial for fixed $k$, this complexity becomes impractical even for moderate values of $n$ and $k$. This limitation is significant because many structural properties in Cops and Robbers only emerge in larger graphs. For example, graphs with both large girth and minimum degree at least three must grow exponentially in size as girth increases, so asymptotic behavior cannot be observed on small instances \cite{Jajcayova2016Improved,bradshaw2020copnumbergraphshigh}. A slow algorithm prevents meaningful exploration of graph properties.

Petr et al. introduced a faster algorithm for determining whether $k$-cops can guarantee capture on a graph \cite{PETR202211}, improving the runtime to $O(kn^{k+2})$. This is achieved through a carefully designed state-space model in which each game-state encodes the positions of all players, as well as a turn indicator. The algorithm applies a reverse breadth-first search (BFS)-like traversal to mark winning states for the cops, starting from base-case copwin states and propagating backward through the state-space tree.

Although Petr et al. provided a detailed theoretical description of their faster $k$-copwin algorithm, they did not present an implementation. Their omission provides an opportunity for other researchers to explore the practical behavior of the algorithm. Reproducing the algorithm is important for research because it helps validate the correctness of the theory and the Cops and Robbers game has many open problems, so developing this tool helps advance research in the field.

\subsection{Contribution}
During the implementation process we clarified the consistency between the theoretical description and practical demands of the implementation. The algorithm models each game configuration as a tuple that encodes the positions of the robber and the cops, along with a turn indicator specifying which player moved last. One of the central data structures in the algorithm is a counter array used to track potential escape options for the robber. The pseudocode initializes the counter array over states where it is the cops' move, rather than the robber's, which is inconsistent with how the state-space is defined earlier in the paper. In order to implement this algorithm, it was necessary to realign this variable, which had other impacts through the rest of the implementation. Our implementation produces results consistent with exhaustive enumeration and theoretical expectations. Implementation, by contrast, demands full precision and reveals ambiguities that may go unnoticed in a purely theoretical setting. This discrepancy could only be noticed by completing a full implementation and speaks to the importance of implementing theoretical algorithms for practical use.

In this paper, we present a Python implementation of the 2022 algorithm by Petr et al., along with a detailed analysis of the differences between the original pseudocode and practical implementation. We provide the full Python code, discuss its implications for both theory and practice, and demonstrate empirically that our implementation restores consistency with expected results on small graphs. By bridging theoretical design with practical implementation, our work reinforces the value of implementation as a means of validating algorithmic research.

\section{Algorithm Overview}

The algorithm introduced by \cite{PETR202211} presents a significant advancement in determining whether a graph is $k$-copwin by representing the  Cops and Robbers as a state-space problem. The algorithm simulates the game of Cops and Robbers on a finite, state graph by encoding each possible game configuration in a state, which is represented as a tuple. Here, they consider each cop $C_1, \ldots C_k$ and the robber $R$ as different player ``pieces'', and each will move in their own turn one at a time, repeating every $k+1$ turns. This is analogous to the cops and robber alternating turns and all cops being allowed to move simultaneously, we simply consider the cop turn to take place over $k$ turns. Each state is represented as a tuple of the form
\begin{center}
    $s=(p_0, p_1, \ldots, p_k, t)$,
\end{center}
where $p_0$ denotes the robber's position, $p_1, \ldots, p_k$ denote the positions of the $k$ cops, and $t$ is a turn indicator. The turn indicator tells us which piece has moved last, so if $t=k$ then it is the robber's turn, and if $t \neq k$ then it is the cop ($t+1$)'s turn. There are three main data structures used: \texttt{COPS\_WIN}, a boolean array indexed by states; \texttt{COUNTER}, an array of non-negative integers indexed by states when it is the robber's turn; and a \texttt{QUEUE} of copwin states.

The main idea of the algorithm is to use reverse BFS through the state space to identify the positions from which the cops can guarantee a win. The algorithm begins by initializing a boolean array, \texttt{COPS\_WIN}, with all states set to 0, indicating that no states are initially marked as winning for the cops. To mark a state as copwin, it is set to 1 in the \texttt{COPS\_WIN} array. The algorithm then identifies the final states where the robber shares a vertex with any of the cops and adds them to the \texttt{QUEUE}. Next, it processes the queued states in the order they were added. For each state $s = (p_0, p_1, \ldots, p_k, t)$, it looks at its predecessors. That is, its incoming neighbours in the state graph. If $t \neq 0$, it means that a cop move reached $s$, so all its predecessors are marked as copwin and added to the queue, unless they have already been marked. Each state is enqueued at most once.

A central data structure used in this process is the \texttt{COUNTER}, which tracks the number of escape options available to the robber in each robber-turn state. It is initialized with the number of out degrees for a given state, that is, $1+deg_G(p_0)$. Here, $G$ represents the input graph, which is an undirected graph with no self-loops that represents state transitions. The 1 is added because the robber can choose to remain in its current position as an available escape option. The \texttt{COUNTER} tracks how many safe moves the robber has from each state $s$. When a state is marked as copwin on the robber's turn, we reduce the \texttt{COUNTER} of each predecessor $s'$ (a state leading to $s$) by 1, since moving to $s$ is no longer safe for the robber. If the \texttt{COUNTER} at $s'$ reaches zero, it means the robber has no safe moves left from $s'$, so $s'$ is also marked as a copwin state.

After all reachable states have been processed, the algorithm checks whether there exists an initial configuration of the $k$ cops (i.e., a fixed $(p_1, \ldots, p_k)$) such that for every starting position of the robber, those game-states are marked as winning. If such a configuration exists, then the graph is declared $k$-copwin. Otherwise, no guaranteed winning strategy exists for the cops, and the graph is not $k$-copwin.

\section{Implementation}


We implemented the algorithm from \cite{PETR202211} in Python, translating the pseudocode into executable code. The algorithm models the Cops and Robbers game as a finite state space, where each state encodes the current positions of all players and a turn indicator. Our implementation computes whether $k$-cops can guarantee the capture of a robber on a given graph by propagating winning states backwards to look at predecessors. It starts from terminal positions where the robber is caught and uses breadth-first traversal to deduce prior game-states that are also winning.

\subsection{Design Choices and Performance Considerations}

Several dictionaries were initialized to store win/loss information (\texttt{COPS\_WIN}), unprocessed move counts for the robber (\texttt{COUNTER}), and state transitions (\texttt{PREDS}). We also maintain a directed graph of state transitions (\texttt{STATE\_GRAPH}) for optional analysis or visualization.

\subsection{Initialization of the Data Structures}
The original algorithm in \cite{PETR202211} compresses the initialization of key data structures into a single line of pseudocode. Based on the descriptions in the paper, we constructed the following data structures in Python:

\begin{lstlisting}[language=Python, caption={Initialization of data structures}]
def is_k_copwin(graph, k):

    n = len(graph)
    if k >= n:
        return True

    COPS_WIN = {}
    COUNTER = {}
    QUEUE = deque()
    vertices = list(range(n))
    PREDS = {}
    STATES = []

    for pos in (product(vertices, repeat=(k + 1))):
        for t in range(0, k+1):
            s = (*pos, t)
            STATES.append(s)
            COPS_WIN[s] = 0
            if t == k:
                COUNTER[s] = len(graph[pos[0]]) + 1
            else:
                COUNTER[s] = 0
\end{lstlisting}

\noindent
The \texttt{PREDS} list is initialized to store predecessor states, which are required during the reverse BFS performed in the while loop explained in Section~\ref{subsec:BFS}. We iterated over all possible configurations of player positions $(p_0, p_1, \ldots, p_k)$ and all possible turn indicators $t \in {0, \dots, k}$, where $t=k$ indicates that it is the robber's turn to move next and the robber's position is at $p_0$. Notice, while the original paper initializes the counter array using $t = 0$, our implementation corrects this to $t = k$ to properly identify robber moves. The robber's turn occurs at $t=k$ because this indicates that the $k^{th}$ cop has moved last. This subtle discrepancy affects the logic of the algorithm, particularly the initialization of \texttt{COUNTER}, since the robber's turn is the only case where unprocessed escape options must be tracked. Correcting this index ensures that the iteration logic correctly handles robber moves and maintains the intended game dynamics throughout state propagation.

\subsection{Identify Final Winning States}

After all data structures are initialized, the winning terminal states, where the robber is already caught, are enqueued for reverse BFS:

\begin{lstlisting}[language=Python, caption={Enqueuing winning terminal states}]
    for s in STATES:
        for i in range(1, k+1):
            robber = s[0]
            cops = s[1:-1]
            t = s[-1]
            if t == i and robber in cops:
                QUEUE.append(s)             
                COPS_WIN[s] = 1
                break
\end{lstlisting}

\noindent
The above code initializes and identifies all final winning states for the cops, that is, the game-states where the robber is already caught. The nested loop iterates over every possible position of the robber and $k$ cops ($n^{k+1}$ combinations), and checks whether the robber occupies the same vertex as any of the cops. If the state has not yet been marked as copwin, then it is added to the QUEUE to begin the reverse BFS process and is marked as copwin.

\subsection{State Transition}
The original pseudocode does not explicitly define how to compute the predecessor states during the reverse BFS phase of the algorithm. In our implementation, we introduce a helper function, predecessors(state), to handle this operation. A game-state is represented as a tuple of the form $s=(p_0, p_1, \ldots, p_k, t)$, where $t$ denotes whose turn it is. The turn indicator, $t$, is the robber ($t = k$) or one of the $k$ cops ($0 \leq t \leq k-1$). Given a current state, the function determines all valid previous states from which the game could have transitioned into it. To do so, it considers every legal position the current player could have come from, which includes all neighbours of their current position, as well as the option to remain in place. For each such prior position, the function constructs a new state by updating the corresponding player's position and rotating the turn index backward using $(t - 1) \bmod (k + 1)$. This ensures proper reversal of the turn order during propagation and enables a breadth-first traversal of the state space from known winning positions.

\begin{lstlisting}[language=Python, caption={Definition of \texttt{predecessors()} function}, label={lst:predecessors}]
def predecessors(state):
    positions = list(state[:-1])
    t = state[-1]
    preds = []

    for prev_pos in graph[positions[t]] + [positions[t]]:
        prev_positions = positions.copy()
        prev_positions[t] = prev_pos
        preds.append(tuple(prev_positions + [(t - 1) % (k + 1)]))

    return preds
\end{lstlisting}

\noindent
Line 2 shows how we extract the positions of all players from the state, discarding the last element $t$, and store them in a list. This list allows us to update the positions as we generate predecessor states. Line 3 retrieves the turn index by identifying which player has just moved. The for loop then considers all positions that the current player could have moved from in the previous step, in order to reach their current position. This includes each neighbour of their current location (\texttt{graph[positions[t]]}), as well as the possibility of remaining in place (\texttt{positions[t]}), as seen in line 6. For each of these possible prior positions, we create a copy of the current position list using \texttt{positions.copy()}. This ensures we do not overwrite the original list, which is used in each iteration of the loop. We then update the copy, \texttt{prev\_positions}, to reflect that the active players had been at this earlier location by setting \texttt{prev\_positions[t] = prev\_pos}. Finally, we append a new state tuple to the list of predecessors. This consists of the modified positions list followed by the rotated turn index, $(t - 1) \% (k + 1)$, which ensures that the game's turn sequence is reversed correctly. The result is a complete list of valid predecessor states that could have led to the current one under legal game transitions.

\subsection{Reverse Breadth-First Search} 
\label{subsec:BFS}
State transitions are computed, and winning states are propagated backwards in a breadth-first manner. The following code shows how the algorithm processes the queue of known winning states and updates the surrounding state space accordingly:

\begin{lstlisting}[language=Python, caption={Reverse BFS propagation of winning states}]
while QUEUE:
    s = QUEUE.popleft()
    _, *_, t = s
    if t != 0:
        for pred in predecessors(s):
            if COPS_WIN[pred] == 0:
                QUEUE.append(pred) 
                COPS_WIN[pred] = 1
                PREDS[pred] = PREDS.get(pred, []) + [s]
    else:
        for pred in predecessors(s):
            COUNTER[pred] -= 1
            if COUNTER[pred] == 0:
                if COPS_WIN[pred] == 0: 
                    QUEUE.append(pred)
                    COPS_WIN[pred] = 1
                    PREDS[pred] = PREDS.get(pred, []) + [s]
\end{lstlisting}

\noindent
This loop forms the core of the BFS process. Each game state $s$ is dequeued from the \texttt{QUEUE} in First In First Out (FIFO) order. All states in the \texttt{QUEUE} will be copwin, see Lemma 2.1 from Petr et al. paper \cite{PETR202211}, so all states dequeued will also be copwin. The unpacking syntax $\_, *\_, t = s$ is used to extract the turn index while ignoring the positional details for brevity.

\textbf{Cop just moved ($t \ne 0$):} If it is a cop’s turn, then any predecessor state from which the game could have transitioned into $s$ is also a winning state for the cops. These predecessor states are identified using the \texttt{predecessors()} function. For each such predecessor:
\begin{itemize}
    \item It is checked if it has not yet been marked as copwin.
    \item It is added to \texttt{QUEUE} for further exploration of its predecessors.
    \item It is marked as a winning state in \texttt{COPS\_WIN}.
    \item The successor relationship is recorded in \texttt{PREDS} for potential backtracking.
\end{itemize}

\textbf{Robber just moved ($t = 0$):} For the robber’s turn, a predecessor is only considered winning for the cops if all of the robber’s escape options lead to states that are already winning for the cops. The \texttt{COUNTER} dictionary tracks how many of these options are still unresolved. For each predecessor:
\begin{itemize}
    \item The counter is decremented to reflect one fewer unexplored option.
    \item If the counter reaches zero and the state is still marked as a loss, it is added to \texttt{QUEUE}, and marked as winning.
\end{itemize}

This propagation continues until the queue is empty, at which point all states from which the cops can force a guaranteed win have been identified.

\subsection{Final Check for k-copwin}

Finally, a global check determines whether every initial robber position can be captured from some cop starting position:

\begin{lstlisting}[language=Python, caption={Final check over initial states}]
for cop_pos in product(vertices, repeat=k):
    if all(COPS_WIN[(robber_pos, *cop_pos, 0)] == 1 for robber_pos in vertices):
        return True
return False
\end{lstlisting}

\noindent
This final loop checks whether $k$ cops can guarantee a win from some initial configuration. The outer loop enumerates all $k$-tuples of starting cop positions via a Cartesian product over the vertex set. For each such tuple, the inner condition tests whether the cops win from every possible initial robber position. The game-state is represented by the tuple $(robber\_pos, *cops\_pos, 0)$, where $robber\_pos$ is the robber's vertex, $*cops\_pos$ are the cops’ positions, and $0$ indicates that the robber has just moved into a position. If such a configuration exists, the graph is $k$-copwin; otherwise, it is not.

\subsection{Full Implementation in Appendix}
For completeness, the full Python implementation of the \texttt{is\_k\_copwin} algorithm, including state initialization, predecessor generation, and reverse BFS logic, is provided in Appendix~\ref{appendix:fullcode}.

\section{Implementation Challenges}
The process of translating the algorithm proposed by Petr et al.\cite{PETR202211} into a functional Python implementation highlighted several areas where theoretical abstraction left room for interpretive choices during implementation. As is often the case with high-level pseudocode, certain components related to state management were intentionally presented in a concise and abstract manner to emphasize the core algorithmic ideas rather than implementation specifics.

One of the primary challenges we encountered involved the lack of explicit initialization for key data structures such as the win-state dictionary, propagation counter, and processing queue. In the original pseudocode, these structures are written in a very condensed way, using just a single line that hides the details of how the full state space is actually handled behind the scenes. As a result, it was not immediately clear how to populate and handle these structures in practice. 

In our implementation, we made these structures explicit. The win-state dictionary \texttt{COPS\_WIN} and the counter array \texttt{COUNTER} are initialized for every possible configuration of robber and cop positions across all turn indices, where $t \in \{0, 1, ..., k\}$. This required iterating over the full Cartesian product of vertex indices to ensure full coverage of the state space. During this process, we also identified what appears to be a minor indexing discrepancy in the pseudocode: the initialization condition uses $t = 0$ to represent the robber’s move, whereas our interpretation of the turn rotation logic suggested this should be $t = k$ instead. This adjustment ensured that the state space was correctly aligned with the intended sequencing of turns.

Another interpretive gap we addressed was the absence of a definition for how predecessor states should be generated during the reverse BFS phase. To fill this gap, we introduced a predecessor function that enumerates all valid previous game-states from which a given state could come. This involved identifying the current player based on the current turn index and generating all legal transitions by considering both adjacent vertices and the possibility of remaining in place. By defining this function, we allowed for consistent traversal of the state space and accurate reverse BFS of winning conditions.

These adjustments show the difference between how ideas are described in theory and how they need to be written in code.  By filling in the missing details, we stayed true to the algorithm’s intent while making the logic executable.

\section{Logical Correctness for Proof}


We base our correctness argument on Lemma 2.2 from the original paper \cite{PETR202211}, which establishes that once the algorithm finishes
\[
\mathrm{COPSWIN}(s) = 1 \quad \Leftrightarrow \quad s \text{ is copwinning}.
\]
The proof of Lemma 2.2 relies on a correct interpretation of the turn index $t$ in a state $s = (p_0, p_1, \dots, p_k, t)$, where $p_0$ is the robber’s position, $(p_1, \dots, p_k)$ are the cops’ positions, and $(t \in \{0, 1, \dots, k\})$ indicates who just moved.

In the original description, the authors implicitly equate $t = 0$ with `robber just moved,' in Section 2 of the article. However, this is inconsistent with the rest of the algorithm and the structure of the state space. In fact, the correct condition for the robber’s turn is $t = k$.

This correction follows from the modular update rule $t' = (t+1) \bmod (k+1)$, where turns advance in the order: $\text{Cop 1} \ (t=0), \ \text{Cop 2} \ (t=1), \ \dots \ , \ \text{Cop $k$} \ (t=k-1), \ \text{Robber} \ (t=k).$ Under this convention, when the robber moves $t = k$, the next turn index becomes $0$, corresponding to Cop 1’s move. This indexing matches the actual move sequence used in our implementation and ensures consistency between state transitions, counter updates, and queue processing. With this correction, the reasoning in Lemma 2.2 remains logically sound.

First direction (\(\mathrm{COPSWIN}(s) = 1 \Rightarrow\) $s$ is copwinning):
When a predecessor state $q$ is processed in the while loop, the turn index $t$ determines whether we are propagating backward from the robber’s move ($t = k$) or a cop’s move ($t \neq k$). The contradiction step in the lemma, which requires that $t = k$ to conclude all robber moves lead to winning states, now follows directly from the corrected indexing.

Second direction ($s$ is copwinning \(\Rightarrow\) \(\mathrm{COPSWIN}(s) = 1\)):
The minimal-move counterexample argument also depends on identifying when the robber moves.  
Using \(t = k\) here ensures that decrementing \(\mathrm{COUNTER}(q)\) after processing all possible robber moves is aligned with the actual game order, so the while-loop’s enqueue step triggers exactly as in the original proof.

Thus, there are no required changes to the proof of Lemma 2.2 because every instance of $t = k$ means the robber's turn. The lemma’s equivalence holds for our corrected algorithm.

\section{Experimental Results}

We evaluated the average runtime performance of two implementations of the $k$-copwin algorithm on all connected graphs with 2 to 9 vertices, totaling 273,192 unique graphs. These graphs were sourced from the online database available at ~\href{https://users.cecs.anu.edu.au/~bdm/data/graphs.html}{https://users.cecs.anu.edu.au/~bdm/data/graphs.html}. The connected graphs can be downloaded from this site in graph6 format, which is the input used for this study. Both algorithms operate under identical rules of the Cops and Robbers game: the game is played on a reflexive graph, and the cops move first.

The first implementation, labeled ``2011 Algorithm'', is a publicly available Python implementation from GitHub, ~\href{https://github.com/Jabbath/Cop-Number}{github.com/Jabbath/Cop-Number} \cite{2011}, based on the algorithm described in Bonato and Nowakowski's 2011 book ``The Game of Cops and Robbers on Graphs" \cite{book}. The second implementation, labeled ``2022 Algorithm'', is our Python implementation of the algorithm proposed by Petr et al. in their 2022 paper ``A faster algorithm for Cops and Robbers" \cite{PETR202211}.

It is a known result in the literature that all graphs with up to 9 vertices have a cop number at most 2. This means that for these graphs, the cop number can only be 1 or 2. Since the smallest graph requiring 3 cops has 10 vertices \cite{minon10}, our evaluation on graphs with at most 9 vertices only needs to verify whether the computed cop number is 1 or greater than 1. In other words, the problem of determining the exact cop number simplifies to distinguishing between cop number 1 (the graph is ``copwin") and cop number 2 (not ``copwin"). This simplification allows us to focus the correctness evaluation on verifying whether the algorithm correctly identifies if a graph is copwin or not, instead of checking for higher cop numbers.

We evaluate the two implementations based on the following criteria. First, we assess correctness by comparing the cop numbers computed by each implementation to ensure consistency and agreement across all graphs. Second, we measure the total computation time each algorithm requires to process all graphs with vertex counts ranging from 2 through 9. This allows us to evaluate and compare the relative scalability and efficiency of the two approaches.

Both algorithms were modified so that output would be written to a csv file for easy comparison. Each csv file contained the following headers:
\begin{itemize}
    \item Graph Number: A number from 1 to 273,192.
    \item Graph Size: Number of nodes ranging from 2-9.
    \item Copwin: True if copwin, False otherwise.
    \item Time: Time in seconds.
    \item Graph6: The graph6 encoding for each graph.
\end{itemize}

When comparing the ``Copwin" columns in both csv files for the two algorithms, we confirmed that all results were consistent. In other words, the cop numbers matched exactly across both implementations, allowing us to proceed confidently with the runtime analysis.

\begin{figure}
\includegraphics[width=1.0\columnwidth]{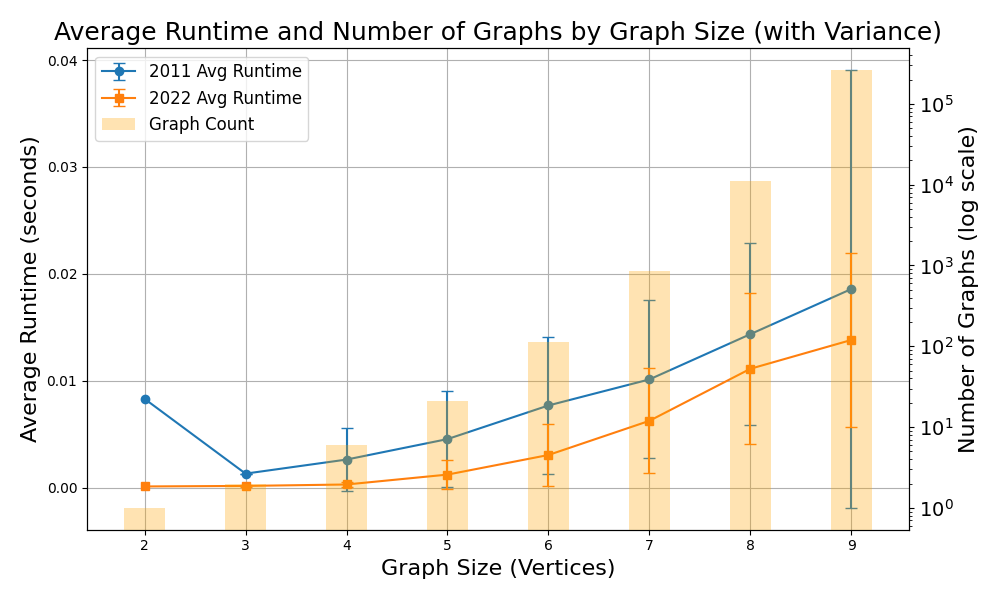}
\caption{{\bf Average runtime comparison with another implementation.} This figure shows the average runtime, in seconds, for both algorithms across all connected graphs of sizes 2-9 vertices. Both algorithms were executed under the same experimental conditions, and the runtimes were averaged over all graphs of a given order. The twin axis displays the number of graphs at each vertex count, providing context for the runtime averages by showing the varying dataset size.}
\label{fig:AverageRuntime}
\end{figure}

Figure \ref{fig:AverageRuntime} presents the average runtime of each algorithm plotted against graph size. A secondary axis displays the corresponding number of graphs for each vertex count. For smaller vertex counts, there are relatively few graphs, while for larger sizes, the number of graphs grows exponentially. Displaying the graph count helps provide context for the runtime measurements, making it clear whether average runtimes are based on a few or many graphs. 

This figure compares the runtime performance of the 2011 and 2022 algorithms on all connected graphs with 2 to 9 vertices. At low vertex counts, the number of connected graphs is very small. At $n=2$ there is 1 graph, at $n=3$ there are 2 graphs, and at $n=4$ there are 7 graphs. The average runtime is computed over fewer instances, so the mean is highly sensitive to the structure of just one or two graphs. The results show that the 2022 algorithm consistently achieves lower average runtimes than the 2011 implementation, confirming its theoretical efficiency gains. The variance of both algorithms generally increases as the graph size increases. Graphs with $n=9$ have the greatest variance for both algorithms, with the exception of the 2011 algorithm at $n=3$. Overall, the 2022 algorithm demonstrates superior scalability on average, though its runtime distribution is less stable across graph instances than that of the 2011 algorithm.

These results confirm that the improvements from the 2022 algorithm translate to practical runtime savings, especially for graph sizes up to $n=9$, where exhaustive enumeration is still feasible.

\section{Implementation Improvements for Efficiency}
\begin{figure}
\includegraphics[width=1.0\columnwidth]{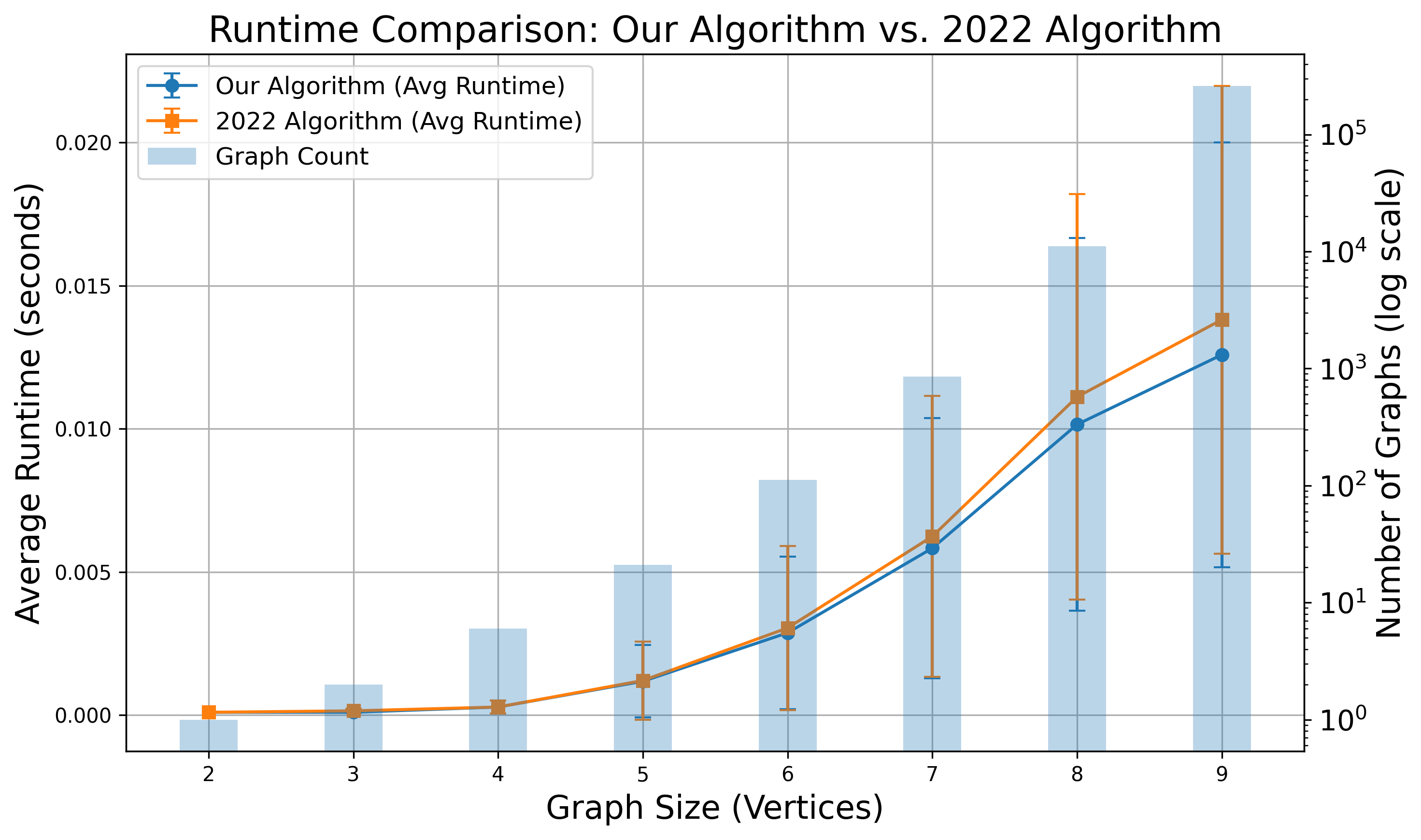}
\caption{{\bf Average runtime comparison with our improved efficiency.} This figure compares the average runtime, in seconds, of the original 2022 algorithm and our optimized implementation across all connected graphs of order 2-9. Runtimes were averaged by graph size under identical conditions. The secondary log-scale axis indicates the number of graphs at each size.}
\label{fig:AverageRuntime2}
\end{figure}

We use an optimized initialization step to mark the winning states in which the robber is already caught. Rather than iterating over the entire state space $(p_0, p_1, ... p_k, t)$, which includes $n^{k+1}$ possible state configurations, we iterate only over the positional configurations $(p_0, p_1, ... p_k)$. For each configuration, we generate only those turn states $(p_0, p_1, ... p_k, t)$ where the robber and the corresponding cop are occupying the same vertex. This avoids states looking at states that cannot be immediate wins. Correctness still holds because both methods mark exactly the states where the robber is already caught. Both versions run in $O(k \cdot n^{k+1})$, but our implementation reduces memory traversal and branch checks, so overall constant factors are reduced. Figure \ref{fig:AverageRuntime2} illustrates the resulting reduction in runtime compared to the original implementation, averaged across all connected graphs of order 2-9.

\section{Discussion}
This project highlights the relationship between theoretical design and practical implementation in research. While the proposed algorithm in the original paper provides a significant theoretical advancement in complexity, the process of translating it into a working Python implementation revealed subtle but important challenges that are common when moving from theory to practice.

The most notable challenge is in the interpretation of pseudocode. Pseudocode in academic papers is often written in high-level abstractions. While this improves readability, it can also obscure important implementation details. In this case, an off-by-one discrepancy in the turn index led to unexpected behaviour until it was resolved. This issue was compounded by hidden assumptions that were not explicitly documented, such as how turn transitions should be handled and how game states should be represented internally. These details are rarely emphasized in theoretical work, but they are essential for a correct, reproducible implementation. By confronting these issues directly, we not only ensured our implementation aligned with the theoretical model but also clarified ambiguities that could hinder replication by others.

This experience also shows that communication between theory and practice matters. Algorithmic research gains value when its results can be adopted, extended, and verified by others. Ambiguities in notation, indexing conventions, or underlying assumptions can create barriers to adoption, and even small semantic differences can lead to diverging results. Our implementation serves as both a validation of the theoretical results and a reference for future researchers, rendering theoretical insights useful in practice.

Finally, we view implementation as a form of validation. When an implementation faithfully reproduces theoretical results, it strengthens confidence in the original algorithm. When discrepancies arise, they can reveal either programming errors or areas where the theory requires further clarification. In this way, implementation serves as a practical tool and also as a test of theoretical soundness. The bridge between theory and practice is bidirectional: theory guides implementation, but implementation also sharpens and strengthens theory.

\section{Conclusions}
In this work, we presented a detailed Python implementation of the 2022 algorithm for deciding whether a graph is $k$-copwin \cite{PETR202211}.  We also proposed improvements to the original algorithm’s description by making implicit assumptions explicit and offering a structured implementation.

Our experimental evaluation compared this implementation to a Python implementation with the same game assumptions \cite{book, 2011} on all connected graphs with 2 to 9 vertices, totalling 273,192 unique graphs. Both algorithms produced identical correctness results across all graphs, confirming that our modifications did not alter the underlying logic. Runtime analysis demonstrated that the 2022 algorithm consistently outperformed the 2011 algorithm, with lower average runtimes for almost every graph size tested and more stable scaling behaviour as $n$ increased.

These findings validate the theoretical efficiency gains of the 2022 algorithm and show that they translate into measurable practical improvements. By making the algorithm executable and evaluating it on an exhaustive dataset, we contribute to its robustness, accessibility, and applicability to larger classes of graphs.

Looking forward, several directions remain open. The current implementation, while functional, could be optimized further to improve runtime and memory efficiency, especially when dealing with large or dense graphs. Additional testing across more complex graph families may yield new insights into the algorithm's scalability and corner cases.

By continuing to bridge theory and practice, we aim to further advance the understanding of pursuit-evasion games and contribute to the growing field of reproducible graph algorithms.

\section*{Funding Declaration}
Erin Meger acknowledges research support from NSERC (2025-05700) and Queen’s University.

\bibliographystyle{ieeetr}
\bibliography{refs}
\appendix
\section{Full Python Implementation}
\label{appendix:fullcode}

The appendix provides the complete Python implementation of our adaptation of the $k$-copwin algorithm, based on the state-space model using reverse BFS. This code includes state initialization, winning-state propagation, and a final decision check for whether $k$ cops can guarantee the capture of the robber on a given graph.

\onecolumn
\begin{lstlisting}[language=Python, caption={Complete implementation of the \texttt{is\_k\_copwin} function}]
from collections import deque
from itertools import product
import networkx as nx

def is_k_copwin(graph, k):
    n = len(graph)
    if k >= n:
        return True

    COPS_WIN = {}
    COUNTER = {}
    QUEUE = deque()
    vertices = list(range(n))
    PREDS = {}
    STATES = []

    for pos in (product(vertices, repeat=(k + 1))):
        for t in range(0, k+1):
            s = (*pos, t)
            STATES.append(s)
            COPS_WIN[s] = 0
            if t == k:
                COUNTER[s] = len(graph[pos[0]]) + 1
            else:
                COUNTER[s] = 0

    for s in STATES:
        for i in range(1, k+1):
            robber = s[0]
            cops = s[1:-1]
            t = s[-1]
            if t == i and robber in cops:
                QUEUE.append(s)             
                COPS_WIN[s] = 1
                break

    def predecessors(state):
        positions = list(state[:-1])
        t = state[-1]
        preds = []

        for prev_pos in graph[positions[t]] + [positions[t]]:
            prev_positions = positions.copy()
            prev_positions[t] = prev_pos
            preds.append(tuple(prev_positions + [(t - 1) % (k + 1)]))

        return preds
        
    while QUEUE:
        s = QUEUE.popleft()
        _, *_, t = s

        if t != 0:
            for pred in predecessors(s):
                if COPS_WIN[pred] == 0:
                    QUEUE.append(pred) 
                    COPS_WIN[pred] = 1
                    PREDS[pred] = PREDS.get(pred, []) + [s]
        else:
            for pred in predecessors(s):
                COUNTER[pred] -= 1
                if COUNTER[pred] == 0:
                    if COPS_WIN[pred] == 0:                        
                        QUEUE.append(pred)
                        COPS_WIN[pred] = 1
                        PREDS[pred] = PREDS.get(pred, []) + [s]

    for cop_pos in (product(vertices, repeat=k)):
        if all(COPS_WIN[(robber_pos, *cop_pos, 0)] == 1 for robber_pos in vertices):
            return True
    return False
\end{lstlisting}
\twocolumn

\end{document}